\documentclass[psamsfonts]{amsart}

\usepackage{amssymb,amsfonts}
\usepackage[all,arc]{xy}
\usepackage{enumerate}
\usepackage{mathrsfs}
\usepackage{dsfont}
\usepackage{graphicx}

\newtheorem{theorem}{Theorem}[section]
\newtheorem{corollary}[theorem]{Corollary}
\newtheorem{proposition}[theorem]{Proposition}
\newtheorem{lemma}[theorem]{Lemma}

\theoremstyle{definition}
\newtheorem{definition}[theorem]{Definition}

\newtheorem{example}[theorem]{Example}

\newtheorem{remark}[theorem]{Remark}

\DeclareMathOperator{\Id}{Id}

\DeclareMathOperator{\Sub}{Sub}

\DeclareMathOperator{\lcm}{lcm}
\DeclareMathOperator{\comp}{comp}

\DeclareMathOperator{\Skel}{Skel}

\DeclareMathOperator{\height}{height}


\numberwithin{equation}{section}

\begin{document}
\title[Finitely Additive, Modular, and Probability Functions]{Finitely Additive, Modular, and Probability Functions on Pre-semirings}

\author[Peyman Nasehpour]{\bfseries Peyman Nasehpour}

\address{Peyman Nasehpour\\
Department of Engineering Science \\
Golpayegan University of Technology \\
Golpayegan\\
Iran}
\email{nasehpour@gut.ac.ir, nasehpour@gmail.com}

\author[Amir Hossein Parvardi]{\bfseries Amir Hossein Parvardi}

\address{Amir Hossein Parvardi\\
Department of Mathematics \\
University of British Columbia \\
Vancouver, BC\\
Canada V6T 1Z2}
\email{parvardi@math.ubc.ca}



\begin{abstract}
In this paper, we define finitely additive, probability and modular functions over semiring-like structures. We investigate finitely additive functions with the help of complemented elements of a semiring. We also generalize some classical results in probability theory such as the Law of Total Probability, Bayes' Theorem, the Equality of Parallel Systems, and Poincar\'{e}'s Inclusion-Exclusion Theorem. While we prove that modular functions over a couple of known semirings are almost constant, we show it is possible to define many different modular functions over some semirings such as bottleneck algebras and the semiring $(\Id(D), + ,\cdot)$, where $D$ is a Dedekind domain.
\end{abstract}

\maketitle

\section{Introduction}

Semirings and other semiring-like algebraic structures such as pre-semirings, hemirings and near-rings have many applications in engineering, especially in computer science and of course in applied mathematics \cite{Golan1999}, \cite{GondranMinoux2008}, \cite{HebischWeinert1998} and \cite{Pilz2011}. Finitely additive and modular functions appear in measure theory \cite{Bogachev2007}, probability theory \cite{Ghahramani2005}, lattice and Boolean algebra theory \cite{Graetzer2011}, module theory \cite{AndersonFuller1992}, and a couple of other branches of mathematics as we will mention in different parts of the paper. Since semiring-like algebraic structures are interesting generalizations of distributive lattices and rings, it is quite natural to ask if finitely additive and modular functions can be defined and investigated in ``pre-semiring theory''. On the other hand, probability functions, as special cases for finitely additive functions, play an important role in probability theory. Therefore it seems quite interesting to see if the classical results of probability theory can be stated and proved in this context as well.

In the present paper, we define and investigate finitely additive, modular, and probability functions over semirings and semiring-like algebraic structures. Since different authors give different definitions for semirings and similar algebraic structures, it is crucial to clarify what we mean by these structures. In this paper, by a ``pre-semiring'', we understand an algebraic structure, consisting of a nonempty set $S$ with two operations of addition and multiplication such that the following conditions are satisfied:

\begin{enumerate}
\item $(S,+)$ and $(S,\cdot)$ are commutative semigroups;
\item Distributive law: $a\cdot (b+c) = a \cdot b + a \cdot c$ for all $a,b,c \in S$.
\end{enumerate}

The definition of pre-semirings is borrowed from the book \cite{GondranMinoux2008} and for more on pre-semirings, one may refer to that.

An algebraic structure that is a pre-semiring and possesses an element that is a neutral element for its addition and an absorbing element for its multiplication, is called a hemiring. Usually, the neutral element of a hermring is denoted by $0$. Any hemiring with a multiplicative identity $1\neq 0$ is called a semiring. For more on semirings and hemirings, one may refer to \cite{Golan1999}.

Let us recall that if $B$ is a Boolean algebra, a real-valued function $\mu$ on $B$ is called \emph{finitely additive}, if $\mu(p \vee q) = \mu(p) + \mu(q)$, whenever $p$ and $q$ are disjoint elements of $B$ \cite[Chap. 15]{Halmos1963}. In Definition \ref{finitelyadditived}, we define a function $f$ from a hemiring $S$ into a commutative semigroup $T$ to be finitely additive if $st=0$ for any $s,t\in S$ implies that $f(s+t) = f(s)+f(t)$. The first section of the present paper is devoted to finitely additive functions. Let us recall that an element $s$ of a semiring $S$ is said to be \emph{complemented}, if there exists an element $c_s\in S$ satisfying $sc_s =0$ and $s+c_s = 1$. The mentioned element $c_s\in S$ is called the complement of $s\in S$. One can easily check that if $s\in S$ has a complement, then it is unique. The complement of $s\in S$, if it exists, is denoted by $s^{\bot}$. Note that if $s$ is complemented, then $s^\bot$ is also complemented and $(s^{\bot})^{\bot} = s$. Also, note that if $s,t \in S$ are complemented, the \emph{symmetric difference} of $s$ and $t$ is defined to be $s \triangle t = s^{\bot}t + s t^{\bot}$ \cite[Chap. 5]{Golan1999}. Finally a semiring $S$ is called zerosumfree, if $s+t = 0$ implies that $s=t=0$ for all $s,t \in S$.

In the first section, we investigate finitely additive functions over complemented elements of semirings and show that if $S$ is a zerosumfree semiring, $T$ is a ring, and $f: S \longrightarrow T$ is a finitely additive and normalized function (i.e., $f(1) = 1$), then for complemented elements $s,t \in S$, the following statements hold (See Proposition \ref{complemented1}, Proposition \ref{normalized1}, and Proposition \ref{normalized2}):

\begin{enumerate}

\item $f(t) = f(ts) + f(ts^{\bot})$.

\item $f(s \triangle t) + 2 f(st) = f(s) + f(t) $, where $s \triangle t = s^{\bot}t + s t^{\bot}$.

\item $f(s) + f(s^{\bot}) = 1$.

\item $f(s^\bot t^\bot) = 1 - f(s) - f(t) + f(st)$.

\end{enumerate}

Similar to the concept of independent events in probability theory, we define some elements $s_1, s_2, \ldots, s_n$ of a pre-semiring $S$ to be independent, if $$f(\prod_{x\in X} x) = \prod_{x\in X} f(x),$$ for any nonempty subset $X$ of $\{s_1,s_2, \ldots, s_n \}$, where $f$ is a function from $S$ into a commutative semigroup $T$ (See Definition \ref{independentd}). After that, in Theorem \ref{independentt}, we prove that if $S$ is a semiring, $R$ a ring, and $f: S \longrightarrow R$ a finitely additive and normalized function, then the following statements for complemented elements $s_1, s_2, \ldots, s_n \in S$ are equivalent:

\begin{enumerate}

\item The elements $s_1, s_2, \ldots, s_n \in S$ are independent,

\item The $2^n$ sets of elements of the form $t_1, t_2, \ldots, t_n$ with $t_i = s_i$ or $t_i = s^{\bot}_i$ are independent,

\item The $2^n$ equalities $f(t_1 t_2 \cdots t_n) = f(t_1) f(t_2) \cdots f(t_n)$ with $t_i = s_i$ or $t_i = s^{\bot}_i$ hold.

\end{enumerate}

In the second section of the present paper, we define probability functions over semirings (Definition \ref{probabilityfunctionsdef}), inspired from the definition of probability functions in probability theory \cite[Chap. I]{Kolmogorov1956}. For defining probability functions, we put an order on the co-domains of those functions and define them as follows:

Let $S$ be a semiring and $T$ an ordered semiring. We define a function $p : S \longrightarrow T$ to be a probability function, if the following properties are satisfied:

\begin{enumerate}

\item $p(s) \geq 0$, for any $s\in S$.

\item $p(1) = 1$.

\item If $s,t \in S$ and $st = 0$, then $p(s+t) = p(s) + p(t)$.

\end{enumerate}

In this section, we prove a couple of nice results for probability functions over semirings similar to probability theory and generalize some of the classical results of probability theory such as Boole's Inequality (Theorem \ref{booleinequality}), the Law of Total Probability (Proposition \ref{totalprobability}), and Bayes' Theorem (Theorem \ref{Bayestheorem}).

Modular functions appear in different branches of mathematics. A modular function is usually a real-valued function $m$ over some objects that two operations $`` + ''$ and $`` \cdot ''$ can apply over those objects and if $s,t$ are of those objects, we have the following ``modular'' equality: $$m(s+t) + m(st) = m(s) + m(t).$$ For example, if $L$ is a lattice of finite length, then $L$ is a modular lattice if and only if the function $\height: L \longrightarrow \mathbb N_0$ is modular, i.e., $$\height(a\vee b) + \height(a\wedge b) = \height(a) + \height(b),$$ for all $a,b \in L$ \cite[Corollary 376]{Graetzer2011}. The third section of the present paper is devoted to modular functions. In Definition \ref{modularfunctionsdef}, we define modular functions over pre-semirings similarly, and examine those functions over different pre-semirings such as Tropical, Max-plus, Truncation, Arctic, and $B(n,i)$ algebras / semirings and show that modular functions over these semirings are almost constant, i.e., they are constant functions over their domains except a finite number of elements in them (See the results from Proposition \ref{ArcticSemiring} to Theorem \ref{semifieldparvardi}).

In different branches of mathematics, there are many interesting examples of modular functions that are not almost constant. We have gathered a couple of them in Example \ref{modularexamples}.

However, we find two interesting examples for nontrivial modular functions. One is given in Proposition \ref{gcdlcm} and the other is given in Theorem \ref{DedekindNP} and is as follows:

If $(G,+)$ is an Abelian group and $D$ is a Dedekind domain, a modular function $f: \Id(D) \longrightarrow G$ can be characterized by the values of $f$ at $(0)$, $D$, and all maximal ideals of $D$.

Finally, we show in Theorem \ref{Poincare} that Poincar\'{e}'s inclusion-exclusion theorem in probability theory holds for arbitrary modular functions over pre-semirings. More precisely, we prove that if $S$ is a multiplicatively idempotent pre-semiring, $T$ a commutative semigroup, and $m: S \longrightarrow T$ is a modular function, then for any elements $s_1, s_2,\ldots, s_n \in S$, the following equality holds: $$m\left(\sum^n_{i=1} s_i\right) + \sum_{i<j} m(s_i s_j) + \sum_{i<j<k<l} m(s_i s_j s_k s_l) + \cdots =$$ $$\sum^n_{i=1} m(s_i) + \sum_{i<j<k} m(s_i s_j s_k)+ \cdots \qquad \text{(Poincar\'{e}'s Formula)},$$ where the number of multiplicative factors in the sums of both sides of the equality is at most $n$.

As a result of Poincar\'{e}'s formula, in Proposition \ref{Poincare2}, we prove that if $S$ is a multiplicatively idempotent pre-semiring, $R$ a ring, and $f: S \longrightarrow R$ a modular function, then for any $n>1$, if $s_1,s_2,\ldots,s_n \in S$ are independent, then $s_1 + s_2 + \cdots + s_{n-1}$ and $s_n$ are also independent.

In the final part of the paper, we also prove that over complemented elements of a semiring, finitely additive functions are modular. In fact, we show that if $S$ is a zerosumfree semiring such that $1+1$ is a complemented element of $S$, $T$ is a commutative semigroup, and $f: S \longrightarrow T$ a finitely additive function, then $$f\mid_{\comp(S)} : \comp(S) \longrightarrow T$$ is modular, where by $\comp(S)$, we mean the set of complemented elements of the semiring $S$ (See Theorem \ref{modular}).

We emphasize that in this paper, the multiplication in all rings, semirings, hemirings, and pre-semirings are commutative, unless otherwise is stated.

\section{Finitely Additive Functions on Hemirings}

Let us recall that if $X$ is a fixed set, then an algebra of sets $\mathcal{A}$ is a class of subsets of $X$ such that $\emptyset, X \in \mathcal{A}$ and if $A, B \in \mathcal{A}$, then $A \cup B \in \mathcal{A}$, $A \cap B \in \mathcal{A}$, and $ A^{\prime} \in \mathcal{A}$. In measure theory, a real-valued function $\mu$ on an algebra $\mathcal{A}$ is called a finitely additive function if for disjoint subsets $A$ and $B$ of $X$, we have $\mu (A \cup B) = \mu(A) + \mu(B)$ \cite{Bogachev2007}. From this definition, we are inspired to define finitely additive functions over hemirings as follows:

\begin{definition}

\label{finitelyadditived}

Let $S$ be a hemiring and $s,t, s_1, s_2,\ldots, s_n \in S$.

\begin{enumerate}

\item The two elements $s$ and $t$ are said to be disjoint, if $st = 0$.

\item The $n$ elements $s_1, s_2, \ldots, s_n$ ($n>1$) are said to be mutually disjoint, if $s_i s_j =0$ for all $1 \leq i < j \leq n$.

\item Let $T$ be a commutative semigroup. We define a function $f : S \longrightarrow T$ to be finitely additive if for any disjoint elements $s$ and $t$ of $S$, we have $$f(s+t) = f(s) + f(t).$$

\end{enumerate}

\end{definition}

Let us mention that there are plenty of examples for finitely additive functions on various algebraic structures defined in different branches of mathematics. Since a modular function $m$ defined over a hemiring with $m(0)=0$ is also a finitely additive function (refer to Theorem \ref{modular}), those examples given in Example \ref{modularexamples} are finitely additive functions as well. Probability functions are also finitely additive. Therefore, the functions, given in Example \ref{probabilityexamples}, are also interesting finitely additive functions. By the way, we give a couple of interesting examples for finitely additive functions over hemirings in the following:

\begin{example}

\label{Finitelyadditiveexample}

Examples of finitely additive functions over hemirings:

\begin{enumerate}

\item Let $\mathcal{I}$ be the collection of finite unions of all intervals of the form $[c,d]$, $[c,d)$, $(c,d]$, and $(c,d)$ where $c,d \in \mathbb R$. It is easy to see that $(\mathcal{I}, \cup, \cap, \emptyset)$ is a hemiring. For any interval $J$ of the form $[c,d]$, $[c,d)$, $(c,d]$, and $(c,d)$, set $\mu(J) = |d-c|$ and for a finite union of disjoint intervals $J_1, J_2,\ldots, J_n$, set $\mu(\bigcup_{i=1}^n J_i) = \sum_{i=1}^n \mu(J_i)$. Then one can easily check that $\mu$ is a finitely additive function on the hemiring $\mathcal{I}$ \cite[Chap. 1]{Bogachev2007}.

\item Let $X$ be an infinite set and $\mathcal{A}$ be the algebra of sets $B \subseteq X$ such that either $B$ or $B^c$ is finite (in some references, when $B^c$ is finite, $B$ is called co-finite). For finite $B$, define $f(B) = 0$, and for co-finite $B$, define $f(B) = 1$. It is easy to see that $(\mathcal{A}, \cup, \cap, \emptyset, X)$ is a Boolean algebra and therefore a hemiring and $f$ a finitely additive function \cite[Problem 2.13]{Billingsley1995}.

 \item Let $T$ be a commutative semigroup and $E$ an entire hemiring, that is $ab=0$ implies that either $a=0$ or $b=0$ for any $a,b \in E$. Any function $f: E \longrightarrow T$ with the property that $f(0) = 0$ is finitely additive and the reason is as follows: If $a,b\in E$ and $ab = 0$, then either $a=0$ or $b=0$ and in each case, $f(a+b) = f(a) + f(b)$.
\end{enumerate}

\end{example}

\begin{proposition}

\label{finitelyadditivep}

Let $S$ be a hemiring, $T$ a commutative semigroup, and $f : S \longrightarrow T$ finitely additive. Then for any disjoint elements $s_1, s_2, \ldots, s_n \in S$, we have $$f(s_1 + s_2 + \cdots + s_n) = f(s_1) + f(s_2) + \cdots + f(s_n).$$

\begin{proof}

By Definition \ref{finitelyadditived}, the statement holds for $n=2$. Now let $n>2$ and suppose that $s_1, \ldots, s_n$ are elements of $S$ such that $s_i s_j =0$ for all $1 \leq i < j \leq n$. It is clear that $$(s_1 + \cdots + s_{n-1}) \cdot s_n =0$$ and therefore $$f((s_1 + \cdots + s_{n-1}) + s_n) = f(s_1 + \cdots + s_{n-1}) + f(s_n).$$ But by induction's hypothesis, $$f(s_1 + \cdots + s_{n-1}) = f(s_1) + \cdots + f(s_{n-1}),$$ and the proof is complete.
\end{proof}

\end{proposition}

Example \ref{Finitelyadditiveexample} shows that a finitely additive function $f : S \longrightarrow T$, where $S$ is a hemiring, will be more interesting if $S$ is not an entire hemiring, i.e., a hemiring with nontrivial zero-divisors. Semirings with nontrivial complemented elements are interesting examples of semirings with nontrivial zero-divisors and suitable for our purpose in this paper. Let us recall that an element $s$ of a semiring $S$ is said to be complemented if there exists an element $c_s\in S$ satisfying $s\cdot c_s =0$ and $s+c_s = 1$. The mentioned element $c_s\in S$ is called the complement of $s$. One can easily check that if $s\in S$ has a complement, then it is unique. The complement of $s\in S$, if it exists, is denoted by $s^{\bot}$. Also, note that the set of all complemented elements of a semiring $S$ is denoted by $\comp(S)$. The set $\comp(S)$ is always nonempty, since $0\in \comp(S)$. It is clear that if $s\in \comp(S)$, then $s^\bot \in \comp(S)$ and $(s^{\bot})^{\bot} = s$. Finally if $s,t \in S$ are complemented, the symmetric difference of $s$ and $t$ is defined to be $s \triangle t = s^{\bot}t + s t^{\bot}$ \cite[Chap. 5]{Golan1999}.

\begin{proposition}

\label{complemented1}

Let $S$ be a semiring, $T$ a commutative semigroup, and $f: S \longrightarrow T$ a finitely additive function. If $s\in S$ is complemented, then the following statements hold:

\begin{enumerate}

\item $f(t) = f(ts) + f(ts^{\bot})$, for any $t\in S$.

\item If, in addition, $t\in S$ is also complemented, we have $$f(s \triangle t) + 2 f(st) = f(s) + f(t). $$

\end{enumerate}

\begin{proof}

$(1)$: Since $(ts)(ts^{\bot})=0$, we have $$f(t) = f(t(s+s^{\bot})) = f(ts + ts^{\bot}) = f(ts) + f(ts^{\bot}).$$

$(2)$: Since $(s^{\bot}t) (s t^{\bot}) = 0$, we have the following: $$f(s \triangle t) + 2 f(st) = f(s^{\bot}t + s t^{\bot}) + 2 f(st) = f(s^{\bot}t) + f(st) + f(s t^{\bot}) + f(st) = f(t) + f(s),$$ and the proof is complete.
\end{proof}

\end{proposition}

A normalized measure in measure theory is a measure $\mu$ with $\mu(1) = 1$ \cite[p. 65]{Halmos1963}. Now it is natural to give the following definition. In current and the next section, we are especially interested in finitely additive normalized functions.

\begin{definition}

Let $S$ and $R$ be semirings. We define a function $f: S \longrightarrow R$ to be normalized if $f(1_S) = 1_R$.

\end{definition}

\begin{proposition}

\label{normalized1}

Let $S$ and $T$ be semirings and $f: S \longrightarrow T$ a finitely additive function such that $f(0) = 0$. Then the following statements are equivalent:

\begin{enumerate}

\item The function $f$ is normalized.

\item For any complemented element $s\in S$, $f(s) + f(s^{\bot}) = 1$.

\end{enumerate}

\begin{proof}

$(1) \Rightarrow (2)$: Since $s\in S$ is complemented, $s s^{\bot} = 0$. Obviously this implies that $f(s) + f(s^{\bot}) = f(s+s^{\bot}) = f(1) = 1$.

$(2) \Rightarrow (1)$: Since $0\in \comp(S)$, by assumption $f(0) + f(0^{\bot}) = 1$. But $f(0) = 0$ and $0^{\bot} = 1$. Therefore, $f(1) = 1$.
\end{proof}

\end{proposition}

Note that in Proposition \ref{normalized1}, it is possible to assume that $T$ is a ring instead of assuming that $f(0) = 0$, since cancellation law of addition in rings implies that $f(0) = 0$ automatically (Refer to Remark \ref{modularfunctions}). Let us recall that the operation $\sqcup$ is defined as $s \sqcup t = s + s^{\bot} t$, where $s\in \comp(S)$ and $t\in S$ \cite[Chap. 5]{Golan1999}.

\begin{proposition}

\label{normalized2}

Let $S$ be a zerosumfree semiring, $R$ a ring, and $f: S \longrightarrow R$ a finitely additive and normalized function. If $s$ and $t$ are complemented elements of $S$, then $$f(s^\bot t^\bot) = 1 - f(s) - f(t) + f(st).$$

\begin{proof}

Let $S$ be a zerosumfree semiring and $s,t \in \comp(S)$. According to the proof of Proposition 5.6 in \cite{Golan1999}, we have $s^\bot t^\bot = (s \sqcup t)^\bot$ and $s \sqcup t \in \comp(S)$. So by Proposition \ref{complemented1}, we have: $f(s^\bot t^\bot) = 1 - f(s \sqcup t) = 1 - f(s + s^{\bot} t) = 1 - f(s) - f(s^{\bot} t) = 1 - f(s) - f(t) + f(st)$.
\end{proof}

\end{proposition}

Now we give the following definition inspired from the concept of independent events in probability theory \cite[p. 115]{Ghahramani2005}.

\begin{definition}

\label{independentd}

Let $S$ and $R$ be pre-semirings and $f: S \longrightarrow R$ a function. We define $s_1, s_2,\ldots, s_n \in S$ to be independent, if $$f(\prod_{x\in X} x) = \prod_{x\in X} f(x),$$ for any nonempty subset $X$ of $\{s_1, s_2, \ldots, s_n \}$.

\end{definition}

\begin{theorem}

\label{independentt}

Let $S$ be a semiring, $R$ a ring, and $f: S \longrightarrow R$ a finitely additive and normalized function. Then the following statements for complemented elements $s_1, s_2, \ldots, s_n$ of $S$ are equivalent:

\begin{enumerate}

\item The elements $s_1, s_2, \ldots, s_n$ are independent,

\item The $2^n$ sets of elements of the form $t_1, t_2, \ldots, t_n$ with $t_i = s_i$ or $t_i = s^{\bot}_i$ are independent,

\item The $2^n$ equalities $$f(t_1 t_2 \cdots t_n) = f(t_1) f(t_2) \cdots f(t_n)$$ with $t_i = s_i$ or $t_i = s^{\bot}_i$ hold.

\end{enumerate}

\begin{proof}

$(1) \Rightarrow (2)$: Let $\{t_1, \ldots, t_m\} \subseteq \{s_2 , \ldots, s_n \}$. By Proposition \ref{complemented1}, we have $$f(t_1 \cdots t_m) = f(t_1 \cdots t_m s_1) + f(t_1 \cdots t_m s^{\bot}_1).$$ But $s_1 , t_1, \ldots, t_m$ are independent. So we have $$f(t_1) \cdots f(t_m) = f(s_1)f(t_1) \cdots f(t_m) + f(t_1 \cdots t_m s^{\bot}_1).$$ This implies that $$f(t_1 \cdots t_m s^{\bot}_1) = f(t_1) \cdots f(t_m) - f(s_1)f(t_1) \cdots f(t_m)=$$ $$(1 - f(s_1)) f(t_1) \cdots f(t_m) = f(s^{\bot}_1)f(t_1) \cdots f(t_m).$$

This means that $s^{\bot}_1$ and $s_2, \ldots, s_n$ are independent. Now it is clear that by mathematical induction, we can prove that $s^{\bot}_1, \ldots, s^{\bot}_i, s_{i+1}, \ldots, s_n$ are also independent.

$(2) \Rightarrow (3)$: Obvious.

$(3) \Rightarrow (1)$: Note that since $f(s) + f(s^{\bot}) = 1$ for any complemented element of $s\in S$, $$f(s_1) f(s_2) \cdots f(s_m) = f(s_1) f(s_2) \cdots f(s_m) (f(s_{m+1}) + f(s^{\bot}_{m+1})) \cdots (f(s_n) + f(s^{\bot}_n)).$$ From the $2^n$ equalities $f(t_1 \cdots t_n) = f(t_1) \cdots f(t_n)$ with $t_i = s_i$ or $t_i = s^{\bot}_i$, we get that $$f(s_1) f(s_2) \cdots f(s_m)=$$ $$f(s_1 \cdots s_m s_{m+1} \cdots s_n) + f(s_1 \cdots s_m s^{\bot}_{m+1} \cdots s_n) + \cdots + f(s_1 \cdots s_m s^{\bot}_{m+1} \cdots s^{\bot}_n).$$ But the elements $s_1 \cdots s_m s_{m+1} \cdots s_n$, $s_1 \cdots s_m s^{\bot}_{m+1} \cdots s_n$, $\ldots$ and $s_1 \cdots s_m s^{\bot}_{m+1} \cdots s^{\bot}_n$ are mutually disjoint, so by using Proposition \ref{finitelyadditivep}, we have: $$f(s_1) f(s_2) \cdots f(s_m) =$$ $$f(s_1 \cdots s_m s_{m+1} \cdots s_n + s_1 \cdots s_m s^{\bot}_{m+1} \cdots s_n + \cdots + s_1 \cdots s_m s^{\bot}_{m+1} \cdots s^{\bot}_n)= $$ $$f(s_1 \cdots s_m (s_{m+1} + s^{\bot}_{m+1}) \cdots (s_n + s^{\bot}_n))=f(s_1 \cdots s_m).$$ The thing we wished to show.
\end{proof}

\end{theorem}

Let us recall that a function $d: X \times X \longrightarrow \mathbb R_{+} \cup \{+ \infty \}$ is defined to be semi-metric, if it satisfies all properties of a metric function except the requirement that $d(x,y) = 0$ implies $x=y$. In other words, it is allowed that distinct points may have a zero distance \cite[Definition 1.4.4.]{BuragoBuragoIvanov2001}. We give a similar definition for $G$-semi-metric functions, where $G$ is an ordered Abelian group.

\begin{definition}

Let $X$ be an arbitrary set and $(G,+,\leq)$ be an ordered Abelian group.

\begin{enumerate}

\item A function $f: X \longrightarrow G$ is defined to be non-negative if $f(x) \geq 0$ for any $x \in X$. If $(X,+,0)$ is a monoid, then the function $f$ is defined to be positive if $f$ is non-negative and $f(x)= 0$ if and only if $x=0$, for any $x\in X$.

\item A function $d: X \times X \longrightarrow G$ defines a $G$-semi-metric on $X$, if $d$ has the following properties:

\begin{enumerate}

\item Non-negativeness: $d(x,y) \geq 0$ and $d(x,x)=0$, for any $x,y \in X$.

\item Symmetry: $d(x,y) = d(y,x)$, for any $x,y \in X$.

\item Triangle inequality: $d(x,z) \leq d(x,y) + d(y,z)$, for any $x,y,z \in X$.

\end{enumerate}

Also, $d: X \times X \longrightarrow G$ defines a $G$-metric on $X$, if $d$ defines a $G$-semi-metric on $X$ and $d(x,y) = 0$ implies that $x=y$, for any $x,y \in X$ (Positiveness).

\end{enumerate}

\end{definition}

\begin{theorem}

Let $S$ be a zerosumfree semiring, $G$ an ordered Abelian group, and $f: S \longrightarrow G$ a non-negative and finitely additive function. Then the function $d(s,t) = f(s \triangle t)$ defines a $G$-semi-metric on $\comp(S)$. Moreover, if the function $f: S \longrightarrow G$ is positive, then $d$ defines a $G$-metric on $\comp(S)$.

\begin{proof}

Let $s,t,u \in \comp(S)$. It is clear that $d(s,t) \geq 0$ and $d(s,s)= f(s \triangle s) = f(0)=0$. Also, the symmetric property can be shown easily: $d(s,t)= f(s \triangle t) = f(t \triangle s) = d(t,s)$. Now we show the triangle inequality holds, that is, $d(s,t) + d(t,u) \geq d(s,u)$.

$d(s,t) + d(t,u) - d(s,u) = f(s \triangle t) + f(t \triangle u) - f(s \triangle u) $

$ = f(s^{\bot}t) + f(s t^{\bot}) + f(t^{\bot}u) + f(t u^{\bot}) - f(s^{\bot}u) - f(s u^{\bot})$

$ = f(s^{\bot}tu) + f(s^{\bot}t u^{\bot}) + f(s t^{\bot} u) + f(s t^{\bot} u^{\bot}) + f(s t^{\bot}u) + f(s^{\bot} t^{\bot}u)$

$ + f(st u^{\bot}) + f(s^{\bot} t u^{\bot}) - f(s^{\bot} t u) - f(s^{\bot} t^{\bot} u) - f(s t u^{\bot}) - f(s t^{\bot} u^{\bot})$

$= 2[f(s t^{\bot} u) + f(s^{\bot}t u^{\bot})] \geq 0.$

Finally, let the function $f: S \longrightarrow G$ be positive and assume that $d(s,t) =0$. Then $f(s^{\bot}t) + f(s t^{\bot}) = 0$. From this we get that $f(s^{\bot}t) = f(s t^{\bot}) = 0$, which in turn implies that $s^{\bot}t = s t^{\bot} = 0$.

Now $s= s\cdot 1 = s(t^{\bot} + t) = st^{\bot} + st = st$. In a similar way, $t = ts$ and therefore $s=t$ and this completes the proof.
\end{proof}

\end{theorem}

\section{Probability Functions on Semirings}

Let us recall that if $\mathcal{A}$ is an algebra of sets over a fixed set $X$, by a probability function, it is meant a non-negative, normalized, and finitely additive function \cite[Chap. I]{Kolmogorov1956}. Inspired by this, we define probability functions over semirings similarly. Note that by an ordered semiring, we mean a semiring $(S, +, \cdot)$ with a partial order $\leq$ on $S$ such that the following statements hold:

\begin{enumerate}

\item If $s \leq t$, then $s+u \leq t+u$ for any $s,t,u \in S$;

\item If $s \leq t$ and $0 \leq u$, then $su \leq tu$ for any $s,t,u \in S$.

\end{enumerate}

For more on ordered semirings, one may refer to \cite[Chap. 2]{Golan2003}.

\begin{definition}

\label{probabilityfunctionsdef}

Let $S$ be a semiring and $T$ an ordered semiring. We define a function $p : S \longrightarrow T$ to be a probability function, if the following properties are satisfied:

\begin{enumerate}

\item $p(s) \geq 0$, for any $s\in S$.

\item $p(1) = 1$.

\item If $s,t \in S$ and $st = 0$, then $p(s+t) = p(s) + p(t)$.

\end{enumerate}

\end{definition}

\begin{example}

\label{probabilityexamples} Here are a couple of examples for probability functions:

\begin{enumerate}

\item Let $X$ be a nonempty finite or countable set. It is clear that the Boolean algebra $(P(X), \cup, \cap)$ has a semiring structure. Let $\{p_x\}_{x\in X}$ be a stochastic sequence, that is, $p_x$ are non-negative real numbers with $\sum_{x} p_x = 1$. One may define a function $p : P(X) \longrightarrow \mathbb R$ with $p(A) = \sum_{k\in A} p_k$. It is a well-known fact in probability theory that $p$ defines a probability function on the semiring $P(X)$ \cite{Grigoryan2010}.

\item Let $a$ and $b$ be two real numbers such that $a < b$ and $I_{ab} = [a,b]$ the closed interval of real numbers. A set $Y \subseteq I_{ab}$ is called to be a subinterval of $I_{ab}$ if $Y$ is in one of the forms $[c,d]$, $[c,d)$, $(c,d]$, or $(c,d)$, where $a \leq c \leq d \leq b$. One may define $S$ to be the collection of all possible finite unions of subintervals of $I_{ab}$. It is, then, easy to check that $(S, \cup, \cap)$ is a semiring. Now we define $p: S \longrightarrow \mathbb R$ in the following way:

\begin{enumerate}

\item $ \displaystyle p(I) = \frac{d-c}{b-a}$, if $I$ is one of the subintervals $[c,d]$, $[c,d)$, $(c,d]$, or $(c,d)$ of $I_{ab}$, where $a \leq c \leq d \leq b$.

\item $p(J) = p(I_1) + \cdots + p(I_m)$, where $J = I_1 \cup \cdots \cup I_m$ and $I_1, \ldots, I_m$ are distinct subintervals of $I_{ab}$ such that $I_{\alpha} \cap I_{\beta}$ is either the empty set or a singleton for any $1 \leq \alpha < \beta \leq m$.
\end{enumerate}

It is, then, easy to see that $p$ defines a probability function on the semiring $(S, \cup, \cap)$ \cite{Ghahramani2005}.

\end{enumerate}

\end{example}

\begin{proposition}

\label{prob1}

Let $S$ be a semiring, $T$ an ordered ring, and $p : S \longrightarrow T$ a probability function. If $s\in S$ is complemented, then the following statements hold:

\begin{enumerate}

\item $0 \leq p(s) \leq 1$.

\item $p(ts) \leq p(t)$, for any $t\in S$.

\end{enumerate}

Moreover, if $S$ is a zerosumfree semiring and $s,t\in S$ are both complemented, then we have the following:

\begin{enumerate}

\item[(3)] $p(st) \geq p(s) + p(t) - 1$.

\end{enumerate}

\begin{proof} $ $

(1): $0 \leq p(s) \leq p(s) + p(s^\bot) = 1$.

(2): Since $p(t) = p(ts) + p(ts^{\bot})$, we have that $p(t) - p(ts) = p(ts^{\bot}) \geq 0$, which means that $p(ts) \leq p(t)$.

(3): By Proposition \ref{normalized2}, we have $p(s^\bot t^\bot) = 1 - p(s) - p(t) + p(st)$. Since $p(s^\bot t^\bot) \geq 0$, $p(st) \geq p(s) + p(t) - 1$.
\end{proof}

\end{proposition}

Let us recall that if $s\in S$ is complemented, then $\sqcup$ is defined as $s \sqcup t = s + s^{\bot} t$. Now we prove the following theorem related to complemented elements of a semiring \cite[Chap. 5]{Golan1999}.

\begin{theorem}

\label{zerosumfree}

Let $S$ be a zerosumfree semiring. Then the following statements are equivalent:

\begin{enumerate}

\item If $s,t \in \comp(S)$, then $s+t \in \comp(S)$,

\item $1+1\in \comp(S)$,

\item If $s,t \in \comp(S)$, then $s+t = s \sqcup t$, for any $s,t \in S$,

\item $(\comp(S), +, \cdot)$ is a Boolean algebra.

\end{enumerate}

\begin{proof}
The assertion $(1) \Leftrightarrow (2) \Leftrightarrow (3)$ holds by \cite[Proposition 5.7]{Golan1999}.

$(3) \Rightarrow (4)$: By Corollary 5.9 in \cite{Golan1999}, if $S$ is a zerosumfree semiring, then $(\comp(S), \sqcup, \cdot)$ is a Boolean algebra. But by assumption $s+t = s \sqcup t$, for any $s,t \in \comp(S)$. Therefore $(\comp(S), +, \cdot)$ is a Boolean algebra.

$(4) \Rightarrow (3)$: Straightforward.
\end{proof}

\end{theorem}

\begin{theorem}[Equality of Parallel Systems]

\label{parallelsystems}

Let $S$ be a zerosumfree semiring such that $1+1\in \comp(S)$, $T$ a ring, and $f : S \longrightarrow T$ a finitely additive and normalized function. If $s_1, s_2,\ldots, s_n \in S$ are complemented and independent, then $$f\left(\sum^n_{i=1} s_i\right) = 1 - \prod^{n}_{i=1} (1-f(s_i)).$$

\begin{proof}
Since $S$ is a zerosumfree semiring such that $1+1\in \comp(S)$, by Theorem \ref{zerosumfree}, $(\comp(S), +, \cdot)$ is a Boolean algebra. But $s_1, s_2, \ldots, s_n \in S$ are complemented, so it is clear that the complement of $\sum^n_{i=1} s_i$ is $\prod^{n}_{i=1} s^{\bot}_i$ and we have $$\sum^n_{i=1} s_i + \prod^{n}_{i=1} s^{\bot}_i = 1 .$$ Since $f$ is a finitely additive and normalized function, we have $$f\left(\sum^n_{i=1} s_i\right) + f\left(\prod^{n}_{i=1} s^{\bot}_i\right) = 1.$$ On the other hand, $s_1,s_2, \ldots, s_n \in S$ are independent, so by Theorem \ref{independentt}, $$f\left(\prod^{n}_{i=1} s^{\bot}_i\right) = \prod^{n}_{i=1} f(s^{\bot}_i).$$ But $f(s^{\bot}_i) = 1 - f(s_i)$ and this finishes the proof.
\end{proof}

\end{theorem}

An important technique in set theory, known colloquially as ``disjointification,'' says that if $A_1, A_2, \ldots, A_n$ are arbitrary sets, then the sets defined by $B_k = A_k \setminus \bigcup^{k-1}_{i=1} A_i$ are mutually disjoint and $\bigcup^n_{i=1} B_i = \bigcup^n_{i=1} A_i$ \cite[Proposition 1.24]{Karr1993}. This technique can be generalized in any Boolean algebra:

\begin{lemma}[Disjointification Technique]

\label{Disjointification}

Let $(B,+,\cdot, 0,1, ^{\prime})$ be a Boolean algebra and $a_1, a_2, \ldots, a_n \in B$. Then the elements $b_k = a_k  a^{\prime}_1 \cdots a^{\prime}_{k-1}$ are mutually disjoint and $b_1 + \cdots + b_n = a_1 + \cdots + a_n$.

\begin{proof}

It is straightforward to see that $b_i$s are mutually disjoint. Now we proceed by induction. It is clear that the assertion holds for $n=1$. Let the assertion hold for $n=k$ and we prove it for the case $n=k+1$. Let $a_1, a_2, \ldots, a_{k+1} \in B$ and set $b_l = a_l  a^{\prime}_1 \cdots a^{\prime}_{l-1}$ for any $ 1 \leq l \leq k+1$. By induction's hypothesis,  $b_1 + \cdots + b_k = a_1 + \cdots + a_k$, so $$b_1 + b_2 + \cdots + b_k + b_{k+1} = a_1 + a_2 + \cdots + a_k + a_{k+1}  a^{\prime}_1 \cdots a^{\prime}_k.$$

Now if we set $s= a_1 + a_2 + \cdots + a_k$, we have that $s^{\prime} = a^{\prime}_1 a^{\prime}_2 \cdots a^{\prime}_k$, and finally $$b_1 + b_2 + \cdots + b_k + b_{k+1} = s+a_{k+1}s^{\prime} = (s+a_{k+1}) (s+s^{\prime}) = s + a_{k+1}.$$ So by induction, the proof is complete. \end{proof}

\end{lemma}

\begin{theorem}[Boole's Inequality]

\label{booleinequality}

Let $S$ be a zerosumfree semiring such that $1+1\in \comp(S)$, $T$ an ordered ring, and $p : S \longrightarrow T$ a probability function. If $s_1, s_2, \ldots, s_n \in S$ are complemented, then $$p(s_1 + s_2 + \cdots + s_n) \leq p(s_1) + p(s_2) + \cdots + p(s_n).$$

\begin{proof}
By Theorem \ref{zerosumfree}, $(\comp(S), +, \cdot)$ is a Boolean algebra. Let $s_1, s_2, \ldots, s_n \in \comp(S)$ and define $t_1 = s_1$ and $t_k = s_k s^{\bot}_1 \cdots s^{\bot}_{k-1}$ for any $1<k \leq n$. It is clear that $t_1, t_2, \ldots, t_n$ are mutually disjoint elements of $\comp(S)$ and by Lemma \ref{Disjointification}, $t_1 + t_2 + \cdots + t_n = s_1 + s_2 + \cdots + s_n$. Moreover, by Proposition \ref{prob1}, $p(t_k) \leq p(s_k)$. Therefore, $$p(s_1 + \cdots + s_n) = p(t_1 + \cdots + t_n) = p(t_1) + \cdots + p(t_n) \leq p(s_1) + \cdots + p(s_n),$$ and the proof is complete.
\end{proof}

\end{theorem}

\begin{definition}

Let $S$ be a semiring, $T$ an ordered semiring, and $p : S \longrightarrow T$ a probability function. If for $t\in S$, $p(t)$ is a multiplicatively invertible element of the semiring $T$, the ``conditional probability of $s$ given $t$,'' denoted by $p(s|t)$, is defined to be $p(s|t) = p(st)/p(t) $.

\end{definition}

The proof of the following proposition is straightforward.

\begin{proposition}

\label{independentp}

Let $S$ be a semiring, $T$ an ordered semiring, and $p : S \longrightarrow T$ a probability function. If for the elements $s,t\in S$, $p(s)$ and $p(t)$ are multiplicatively invertible elements of the semiring $T$, then the following statements are equivalent:

\begin{enumerate}

\item $p(s|t) = p(s)$,

\item $p(t|s) = p(t)$,

\item The elements $s$ and $t$ are independent, i.e., $p(st) = p(s)p(t)$.

\end{enumerate}

\end{proposition}

\begin{proposition}
Let $S$ be a semiring, $T$ an ordered semiring, and $p : S \longrightarrow T$ a probability function. If for an element $t\in S$, $p(t)$ is a multiplicatively invertible element of the semiring $T$, then $p_t: S \longrightarrow T$ defined by $p_t(s)=p(s|t)$ is a probability function, i.e., the following statements hold:

\begin{enumerate}

\item $p_t(s) \geq 0$ for any $s\in S$.

\item $p_t(1) = 1$.

\item If $s_1 s_2 =0$, for $s_1,s_2 \in S$ then $p_t(s_1 + s_2) =  p_t(s_1) + p_t(s_2)$.

\end{enumerate}

\begin{proof}
Straightforward.
\end{proof}
\end{proposition}

\begin{proposition}[Law of Total Probability]

\label{totalprobability}

Let $S$ be a semiring, $T$ an ordered semiring, and $p : S \longrightarrow T$ a probability function. If $t_1, t_2,\ldots, t_n$ are elements of $S$ such that $t_i t_j =0$ for any $1 \leq i < j \leq n$, $t_1 + t_2+ \cdots + t_n = 1$, and $p(t_i)$ is a multiplicatively invertible element of the semiring $T$ for any $1 \leq i \leq n$, then $$p(s) = p(s|t_1) p(t_1) + p(s|t_2) p(t_2)+ \cdots + p(s|t_n) p(t_n).$$

\begin{proof}

It is clear that $s = st_1 + st_2 + \cdots + st_n$ and $(st_i)(st_j)=0$ for any $1 \leq i < j \leq n$. So by Proposition \ref{finitelyadditivep}, we have that $p(s) = p(st_1) + p(st_2) \cdots + p(st_n)$. Since $p(st_i) = p(s|t_i) p(t_i)$ for any $1\leq i \leq n$, $$p(s) = p(s|t_1) p(t_1) + \cdots + p(s|t_n) p(t_n),$$ and the proof is complete.
\end{proof}

\end{proposition}

\begin{corollary}[Law of Total Probability]

Let $S$ be a semiring, $T$ an ordered semiring, and $p : S \longrightarrow T$ a probability function. If $t\in S$ is a complemented element of $S$ such that $p(t)$ and $p(t^{\bot})$ are multiplicatively invertible elements of the semiring $T$, then $p(s) = p(s|t) p(t) + p(s|t^{\bot}) p(t^{\bot})$.

\end{corollary}

\begin{proposition}

Let $S$ be a semiring, $T$ an ordered semiring, and $p : S \longrightarrow T$ a probability function. If for the elements $s\in S$, $t\in \comp(S)$, $p(s)$, $p(t)$, and $p(t^{\bot})$ are multiplicatively invertible elements of the semiring $T$, then the following statements are equivalent:

\begin{enumerate}

\item $p(s|t) = p(s|t^{\bot})$,

\item The elements $s$ and $t$ are independent.

\end{enumerate}

\begin{proof}
$(1) \Rightarrow (2)$: By the Law of Total Probability, we have $$p(s) = p(s|t) p(t) + p(s|t^{\bot}) p(t^{\bot}) = $$ $$p(s|t) p(t) + p(s|t) p(t^{\bot}) = p(s|t) (p(t) + p(t^{\bot})) = p(s|t).$$

$(2) \Rightarrow (1)$: Since the elements $s$ and $t$ are independent, by Proposition \ref{independentp}, we have that $p(s) = p(s|t)$ and $p(s) = p(s|t^{\bot})$.
\end{proof}

\end{proposition}

 A purpose of Bayes' formula in probability theory is to compute $P(B|A)$ in terms of $P(A|B)$. We prove a semiring version of Bayes' Theorem for probability functions as follows:

\begin{theorem}[Bayes' Theorem]

\label{Bayestheorem}

Let $S$ be a semiring, $K$ an ordered semifield, and $p : S \longrightarrow K$ a probability function. If $t_1, \ldots, t_n$ are elements of $S$ such that $t_i t_j =0$ for any $1 \leq i < j \leq n$, $t_1 + \cdots + t_n = 1$, $p(t_i) > 0$ for any $1 \leq i \leq n$, and $s\in S$ such that $p(s)>0$, then $$ \displaystyle p(t_k | s) = \frac{p(s|t_k)p(t_k)}{p(s|t_1)p(t_1) + \cdots + p(s|t_n)p(t_n)}.$$

\begin{proof}
By definition, we know that $ \displaystyle p(t_k | s) = \frac{p(t_k s)}{p(s)}$. Also, $p(t_k s) = p(s t_k) = p(s|t_k) p(t_k)$. On the other hand, by the Law of Total Probability (Proposition \ref{totalprobability}), we have $$p(s) = p(s|t_1)p(t_1) + \cdots + p(s|t_n)p(t_n).$$ Therefore, $$p(t_k | s) = \frac{p(s|t_k)p(t_k)}{p(s|t_1)p(t_1) + \cdots + p(s|t_n)p(t_n)},$$ and the proof is complete.
\end{proof}

\end{theorem}

\section{Modular Functions on Pre-semirings}

Modular functions appear in different branches of mathematics. A modular function is usually a real-valued function $m$ over some objects that two operations $``+"$ and $``\cdot"$ have meaning for those objects and if $s$ and $t$ are of those objects, we have the following equality: $$m(s+t) + m(st) = m(s) + m(t).$$ In this section, we define modular functions on pre-semirings and examine modular functions over different pre-semirings mentioned in literature and show that modular functions over these pre-semirings are almost constant, i.e., they are constant over their domains except for a finite number of points in them. We also prove that there are many modular functions over the semiring $(\Id(D),+,\cdot)$, where $D$ is a Dedekind domain. After that, we give a couple of nontrivial examples for modular functions in Remark \ref{modularfunctions} and Example \ref{modularexamples}. Then we prove Poincar\'{e}'s Inclusion-Exclusion Theorem for modular functions over pre-semirings. We also prove that over complemented elements of a semiring, finitely additive functions are modular.

\begin{definition}

\label{modularfunctionsdef}

Let $S$ be a pre-semiring and $T$ a commutative semigroup. We define a function $m : S \longrightarrow T$ to be modular if $m(s+t) + m(st) = m(s) + m(t)$ for all $s,t \in S$.

\end{definition}

\begin{remark}

$ $

\label{modularfunctions}

\begin{enumerate}

\item \label{modularbasicsemirings}If $S$ is a pre-semiring such that $0$ and $1$ are the additive and multiplicative neutral elements of $S$, respectively, $T$ is a commutative and cancellative monoid, and $f: S \longrightarrow T$ is a modular function, then by choosing $t=1$, from the following modular equality $$ f(s+t) + f(st) = f(s) + f(t),$$ we get that $f(s+1) = f(1)$. This means that $f$ is constant on the basic sub-pre-semiring \{1, 1+1, 1+1+1,\ldots\} of $S$.

\item Let $S$ be a hemiring, $T$ a commutative monoid, and $m : S \longrightarrow T$ a function. It is clear that if $m$ is modular and $m(0) = 0$, then $m$ is finitely additive. The question arises when the inverse is correct. Note that if $T$ is a cancellative monoid and $f$ is finitely additive then $f(0) = 0$ and the proof is as follows: since $0\cdot 0 = 0$, we have that $f(0+0) = f(0) + f(0)$ and since $T$ is cancellative, $f(0) = 0$. In Theorem \ref{modular}, we show that in Boolean algebras modularity of a function $f$ is caused by its finitely additive property.

\item Let $S$ be a totally ordered set. It is easy to check that $(S, \oplus, \odot)$ is a pre-semiring (known as bottleneck pre-semiring), where $a \oplus b = \max\{a,b\}$ and $a \odot b = \min\{a,b\}$ for all $a,b \in S$. Now let $(T,+)$ be a commutative semigroup. Obviously for any function $m: S \longrightarrow T$, we have $$m(x \oplus y) + m(x \odot y) = m(x) + m(y).$$

\end{enumerate}

\end{remark}

Now we start examining modular functions over different pre-semirings. Let $\mathbb N_0$ denote the set $\{0,1,2, \ldots,k,k+1, \ldots\}$, i.e., the set of all nonnegative integers. One may define addition and multiplication over $S = \mathbb N_0 \cup \{- \infty \}$ as $``\max"$ and $``+"$ respectively by considering that $ - \infty < n < n+1$ for all $n\in \mathbb N_0$ and $-\infty + s = -\infty$ for all $s\in S$. It is easily checked that $(\mathbb N_0 \cup \{- \infty \}, \max, +, -\infty, 0)$ is a semiring and in some references, it is known as the \emph{Arctic semiring} \cite[p. 179]{DrosteKuichVogler2009}.

\begin{proposition}

\label{ArcticSemiring}

Let $f: \mathbb N_0 \cup \{- \infty \} \longrightarrow G$ be a function from the Arctic semiring $(\mathbb N_0 \cup \{- \infty \}, \max, +, -\infty, 0)$ to the Abelian group $G$. Then $f$ is a modular function if and only if it is a constant function over $\mathbb N_0$.

\begin{proof}
$\Rightarrow$: By modularity of the function $f$, we have: $$f(\max\{x,y\}) + f(x+y) = f(x) + f(y).$$ If we suppose that $x\geq 0$ and $y=0$, we get that $f(x) + f(x) = f(x) + f(0)$ and this means that $f$ is constant over $\mathbb N_0$.

$\Leftarrow$: On the other hand, if $f$ is constant over $\mathbb N_0$, then $f$ is modular over the Arctic semiring.
\end{proof}

\end{proposition}

Let us recall that the semiring $(T_k, \max, \min \{a+b,k\}, - \infty, 0)$, where $1 \leq k$ and $T_k = \{ - \infty, 0, 1, \ldots, k\}$ is called the \emph{Truncation semiring} .

\begin{proposition}

\label{TruncationSemiring}

Let $G$ be an Abelian group. Then $f: T_k \longrightarrow G$ is modular if and only if $f$ is constant over $T_k - \{-\infty \}$, where by $T_k$, we mean the Truncation semiring.

\begin{proof}
$\Rightarrow$: Let $T_k$ be the Truncation semiring, $G$ an Abelian group, and $f: T_k \longrightarrow G$ a modular function. So by definition, we have $$f(\max \{x,y\}) + f(\min \{x+y,k\}) = f(x) + f(y),$$ for all $x,y \in T_k$ and if we let $y=k$, we have the following: $$f(k) + f(\min \{x+k,k\}) = f(x) + f(k).$$ This implies that $f(\min \{x+k,k\}) = f(x).$ But $\min \{x+k,k\} = k$, if $x \geq 0$. So $f$ is constant over $T_k - \{-\infty \}$.

$\Leftarrow$: Straightforward.
\end{proof}

\end{proposition}

Another important family of finite semirings is mentioned in \cite[Example 1.8]{Golan1999} and we bring it here for the convenience of the reader. Let us recall that if $i$ and $n$ are positive integers such that $i<n$, the addition and multiplication of the finite semiring $$\big(B(n,i) = \{0,1,\ldots,n-1 \},\oplus,\odot,0,1\big)$$ are defined as follows:

The addition $\oplus$ is defined as $x\oplus y = x+y$ if $x+y \leq n-1$ and $x\oplus y = l$ if $x+y > n-1$ where $l$ is the unique number satisfying the conditions $i \leq l \leq n-1$ and $l \equiv x+y \pmod{n-i} $ and multiplication $\odot$ is defined similarly. Now we prove the following interesting result:

\begin{proposition}

\label{parvardiproposition}

Let $G$ be an Abelian group and $i<n$ be positive integers. Then a function $f: B(n,i) \longrightarrow G$ is modular if and only if it is constant over $B(n,i)-\{0\}$.

\begin{proof}

$\Rightarrow$: Let $y=1$. It is clear that if $x<n-2$, then $x+y$ and $xy$  are both less than $n-1$. So $x \oplus 1 = x+1$, and $x \odot 1 = x1 =x$.

Putting $y=1$ and $x<n-2$ in the modular equation, we have $$f(x+1) + f(x) = f(x) + f(1),$$ which implies $$f(x+1)=f(1),$$ for $x=0,1,\ldots,n-2$. Thus $f(x) = f(1)$ for $x=1,2,\ldots,n-1$. This means that $f$ is constant over $B(n,i)-\{0\}$.

$\Leftarrow$: Obvious.
\end{proof}

\end{proposition}

\begin{proposition}

\label{constantmodular}

Let $(G,+, \leq)$ be a totally ordered Abelian group and $T$ a commutative and cancellative monoid. Then $(G,\min, +)$ is a pre-semirng and the only modular function, $m: G \longrightarrow T$, over the pre-semiring $G$ is a constant function.

\begin{proof}
It is straightforward to check that $(G,\min, +)$ is a pre-semiring. Note that by assumption, the following equality holds for all $x,y\in G$: $$m(\min\{x,y\}) + m(x+y) = m(x) + m(y).$$ Now take $x\geq 0$. It is, then, clear that $x \geq -x$ and $\min\{x,-x\} = -x$. So by modularity of $m$ over the pre-semiring $G$, we have $m(-x) + m(0) = m(x) + m(-x)$ and therefore $m(x) = m(0)$. On the other hand, if we let $x < 0$, then $\min\{x,0\} = x$ and it is clear that the modularity of $m$ implies that $m(x) + m(x) = m(x) + m(0)$ and again we have $m(x) = m(0)$. This means that for any $x\in G$, we have $m(x) = m(0)$ and $m$ is a constant function and the proof is complete.
\end{proof}

\end{proposition}

\begin{proposition}

\label{constantmodular2}

Let $S$ be a pre-semiring such that $0$ and $1$ are the additive and multiplicative neutral elements of $S$, respectively. If $1$ has an additive inverse, $T$ is a commutative and cancellative monoid, and a function $f: S \longrightarrow T$ is modular, then $f$ is a constant function.

\begin{proof}
Let $x\in S$ be arbitrary and set $y=1$. So by modularity, we have: $$f(x+1) + f(x\cdot 1) = f(x) + f(1).$$ Since $T$ is a cancellative monoid, we have $f(x+1) = f(1)$. Now since $1$ has additive inverse, we can replace $x$ by $x-1$ and we have $f(x) = f(1)$, which means that the function $f$ is constant on $S$ and this finishes the proof.
\end{proof}

\end{proposition}

\begin{corollary}

\label{constantmodular3}

If $R$ is a commutative ring with identity, $T$ a commutative and cancellative monoid, and a function $f: S \longrightarrow T$ is modular, then $f$ is a constant function.
\end{corollary}

\begin{remark}

\label{exconstantmodular2}

Note that if the multiplicative neutral element $1$ of a semiring $R$ has an additive inverse, i.e., $1 + (-1) = 0$, then by distributive law, we have $r + (-1)r = 0$ for any $r\in R$ and this means that $R$ is a ring. Therefore, in order to show that Proposition \ref{constantmodular2} is really a generalization of Corollary \ref{constantmodular3}, it is important to give an example of a pre-semiring $S$ such that $0$ and $1$ are the additive and multiplicative neutral elements of $S$, respectively, and $1$ has an additive inverse, while $S$ is not a semiring. In fact, the interesting example given in \cite[Example 5.3.1]{GondranMinoux2008} fulfills our purpose and inspires us to give the following result:

\end{remark}

\begin{proposition}

Let $(H,+,0,-,\leq)$ be a nontrivial totally ordered Abelian group and $E$ be the set of all closed intervals $[a,b] \subseteq H$ such that $a\leq 0$ and $b\geq 0$. Define $\oplus$ and $\otimes$ on $E$ by $$[a_1,b_1] \oplus [a_2,b_2] = [\min\{a_1,a_2\}, \max\{b_1,b_2\}],$$ $$[a_1,b_1] \otimes [a_2,b_2] = [a_1 + a_2, a_2 + b_2],$$ for all $a_1,a_2 \leq 0$ and $b_1,b_2 \geq 0$. Then the following statements hold:

\begin{enumerate}

\item $(E,\oplus,\otimes)$ is a pre-semiring and $[0,0]$ is the neutral element for both operations $\oplus$ and $\otimes$.

\item The neutral element for multiplication $\otimes$ in $E$ has an additive inverse.

\item The neutral element $[0,0]$  of addition $\oplus$ is not an absorbing element of $E$ and so $E$ is not a semiring.

\item If $G$ is an Abelian group, then $f: E \longrightarrow G$ is modular if and only if $f$ is constant.

\end{enumerate}

\begin{proof}

$ $

$(1)$ and $(2)$: Straightforward.

$(3)$: $[a,b] \otimes [0,0] = [a,b]$ for all $[a,b] \in E$.

$(4)$: This holds by Proposition \ref{constantmodular2}.
\end{proof}

\end{proposition}

\begin{proposition}

\label{Litvinov}

Let $G$ be an Abelian group and $R = (\mathbb R \times \mathbb R) \cup \{- \infty\}$. Define operations of $R$ as follows:
\begin{enumerate}
\item $(a,b) \oplus (a',b') = (\max\{a,a'\}, \max\{b,b'\})$ for all $a,b,a',b' \in \mathbb R$;
\item $(a,b) \odot (a',b') = (a+a',b+b')$ for all $a,b,a',b' \in \mathbb R$;
\item $(-\infty) \oplus r = r \oplus (-\infty)=r$ for all $r \in R$;
\item $(-\infty) \odot r = r \odot (-\infty)=-\infty$ for all $r \in R$.
\end{enumerate}
Then $(R,\oplus,\odot)$ is a semiring and any modular function $f: R \longrightarrow G$ is constant.

\begin{proof}
According to \cite[Example 1.25]{Golan1999}, $R$ is a semiring and it is clear that $ \mathds{1}=(0,0)$ is the identity element of the multiplication of $R$. Let $x,y \in \mathbb R$. By modularity of $f$, we have $$f((x,y) \oplus (0,0))+f((x,y) \odot (0,0))=f(x,y)+f(0,0), $$
which means that $$f(\max\{x,0\}, \max\{y,0\})=f(0,0),$$ and clearly this implies that if $x,y \geq 0$, then $$f(x,y) = (0,0).$$

On the other hand, since $(x,y) \oplus (-x,-y) = (|x|,|y|)$ and $(x,y) \odot (-x,-y) = (0,0)$, modularity of $f$ implies that $$f(|x|,|y|) + f(0,0) = f(x,y) + f(-x,-y).$$ Therefore, if we let $x,y \geq 0$, we have $$f(-x,-y) = f(0,0).$$

Now, let $x,y\geq 0$. Then it is clear that $(x,y) \oplus (x,-y) = (x, y)$ and $(x,y) \odot (x,-y) = (2x,0)$ and again by modularity of $f$, we have $$f(x,y) + f(2x,0) = f(x,y) + f(x, -y).$$ Also, $f(2x,0) = f(0,0)$. So, $f(x,-y) = f(0,0)$. Similarly, one can see that $f(-x,y) = f(0,0)$ and the proof is complete.
\end{proof}

\end{proposition}

\begin{proposition}

Let $(G,+)$ be an Abelian group. Then the following statements for a function $f: \mathbb N_0 \cup \{+ \infty \} \longrightarrow G$ are equivalent:

\begin{enumerate}

\item The function $f: \mathbb N_0 \cup \{+ \infty \} \longrightarrow G$ is a modular function from the tropical algebra $(\mathbb N_0 \cup \{+ \infty \}, \min,+)$ to the Abelian group $G$,

\item The function $f$ is constant over $\mathbb N$.

\end{enumerate}

\begin{proof}
$(1) \Rightarrow (2)$:  Let $y=1$ and $x \in \mathbb N$ and rewrite the modular equality for the tropical algebra:
\begin{equation*}
f(\min\{x,1\}) + f(x+1) = f(x) + f(1).
\end{equation*}
Since $x \in \mathbb N$, $\min\{x,1\}=1$, we have $f(\min\{x,1\})=f(1)$, which gives
\begin{equation*}
f(1) + f(x+1) = f(x) + f(1).
\end{equation*}
Canceling $f(1)$ from both sides, we have $f(x+1)=f(x)$ for all $x \in \mathbb N$, which means $f$ is constant over $\mathbb N$.

$(2) \Rightarrow (1)$: Since $f$ is constant over $ \mathbb N$, we may choose constants $k,c,l \in G$ such that
\begin{equation*}
f(x) =
\begin{cases}
k \quad \mbox{if } x=0\\
c \quad \mbox{if } x\in \mathbb N\\
l \quad \mbox{if } x=+\infty\\
\end{cases}
\end{equation*}
It is straightforward to see that the function $f$ satisfies the modular equality and this finishes the proof.
\end{proof}

\end{proposition}

\begin{theorem}

\label{semifieldparvardi}

Let $G$ be an Abelian group and $(K,+,\cdot, \leq)$ a totally ordered positive semifield with the property that for $y\geq 1$, there is an $x\geq 0$ such that $y=x+1$ for all $x,y \in K$. Then the following statements for the function $f: K \longrightarrow G$ are equivalent:

\begin{enumerate}

\item The function $f$ is modular,

\item The function $f$ is constant over the set of positive elements $X = \{x > 0 : x\in K\}$ of $K$.

\end{enumerate}

\begin{proof}
$(2) \Rightarrow (1)$: Straightforward.

$(1) \Rightarrow (2)$: Since $f$ is a modular function, it satisfies the following equation: $$f(x+y)+ f(xy)=f(x)+f(y).$$
Let $y=1$ and rewrite the functional equation: \[f(x+1) + f(x) = f(x) + f(1),\] which implies that $f(x+1)=f(1)$ for all $x\geq 0$. Since for all $y\geq 1$, there is an $x\geq 0$ such that $y=x+1$, we have that \begin{equation}\label{eq1}
	f(x)=f(1) \quad \forall x \geq 1.
\end{equation} Now replace $\displaystyle y = \frac{1}{x}$ in the original equation to obtain
\begin{equation}\label{eq2}
f\left(x+\frac{1}{x}\right)+f(1)=f(x)+f\left( \frac{1}{x} \right).
\end{equation}
	 Assume $0 < x\leq 1$, then $\displaystyle \frac{1}{x} \geq 1$ and from equation \eqref{eq1}, we have $\displaystyle f\left( \frac{1}{x} \right)=f(1)$. Replacing this in equation \eqref{eq2}, we get $\displaystyle f\left(x+\frac{1}{x}\right)+f(1)=f(x)+f(1)$. Removing $f(1)$ from both sides, we have
\begin{equation*}\label{eq3}
f\left(x+\frac{1}{x}\right) = f(x) \quad \forall~{ } 0 < x \leq 1.
\end{equation*}
	Since $\displaystyle x + \frac{1}{x} \geq 1$, we have $\displaystyle f\left(x+\frac{1}{x}\right)=f(1)$.
	 Therefore, $f(x) = f(1)$ for all $0 < x \leq 1$. Finally, $f(x)=f(1)=c$ for some constant $c\in T$ and all $x\in X$ and the proof is complete.
\end{proof}
\end{theorem}

\begin{example}

\label{semifieldparvardiexamples}

In order to show the importance of Theorem \ref{semifieldparvardi}, let us give a couple of nice examples for semifields satisfying the conditions of Theorem \ref{semifieldparvardi}.

\begin{enumerate}

\item It is clear that $(\mathbb Q^{\geq 0}, +, \cdot, \leq)$ is a totally ordered positive semifield with the property that for $y\geq 1$, there is an $x\geq 0$ such that $y=x+1$ for all $x,y\in \mathbb Q^{\geq 0}$.

\item It is easy to check that if $(K,+,\cdot, \leq)$ is a totally ordered positive semifield, then the semifield $(K, \max, \cdot, \leq)$ is a totally ordered positive semifield with the property that for $y\geq 1$, there is an $x\geq 0$ such that $y=\max\{x,1\}$ for all $x,y \in K$.

\item If $(G,+, \leq)$ is a totally ordered Abelian group, then  $(G \cup \{- \infty\},\max,+, \leq)$, known as $G$-max-plus algebra, is a totally ordered positive semifield with the property that for all $y \geq 0$, there is an $x \geq - \infty$ such that $y = \max \{ x,0 \}$.

\item For nonnegative real numbers $a,b$, define $a \oplus_h b = (a^{1/h} + b^{1/h})^h$, where $h$ is a fixed positive real number. It is easy to see that $S_h = (\mathbb R^{\geq 0}, \oplus_h, \cdot)$ is a totally ordered positive semifield with the property that for $y\geq 1$, there is an $x\geq 0$ such that $y=x \oplus_h 1$ for all $x,y$. For more on the semiring $S_h$, one can refer to \cite{Maslov1986}.

\end{enumerate}

\end{example}

\begin{remark}

Let $\mathbb N$ be the set of natural numbers (positive integers) and $a,b\in \mathbb N$. Then it is easy to check that if $\lcm\{a,b\}$ and $\gcd\{a,b\}$ denote the least common multiple and greatest common divisor of the natural numbers $a$ and $b$, then $(\mathbb N, \lcm, \gcd)$ is a pre-semiring \cite[Example 6.5.8]{GondranMinoux2008} and $\lcm\{a,b\} \cdot \gcd\{a,b\} = a \cdot b$ \cite[Theorem 6.4]{Loo-Keng1982}. This leads us to define a nontrivial modular function as follows:

\end{remark}

\begin{proposition}

\label{gcdlcm}

Define operations $\oplus$ and $\odot$ on $\mathbb N$ as $a \oplus b = \lcm\{a,b\}$ and $a \odot b = \gcd\{a,b\}$. Then the function $f: \mathbb N \longrightarrow \mathbb R$ defined as $f(x) = \log x$ is modular, i.e., $$f(a \oplus b) + f(a \odot b) = f(a) + f(b).$$

\end{proposition}

Let us recall that a semiring $S$ is called simple if $s+1 = 1$ for any $s\in S$ \cite[p. 4]{Golan1999}. We first prove the following lemma for modular functions over simple semirings. After that, we are able to give another interesting nontrivial modular function.

\begin{lemma}

\label{simplesemiring}

Let $S$ be a simple semiring, $G$ an Abelian group, and $f: S \longrightarrow G$ a modular function. Then $f(x^m y^n) = f(xy)$ for any $x,y \in S$ and $m,n \in \mathbb N$.

\begin{proof}
Let $x,y \in S$ and apply the modular equality for the elements $y, yx$. So, $$f(y+yx) + f(y \cdot yx) = f(y) + f(yx).$$ Also, $y + yx = y(x+1) = y$. Therefore, $f(xy^2) = f(xy)$. Now by induction on $m$ and $n$, it is easy to prove that $f(x^m y^n) = f(xy)$.
\end{proof}

\end{lemma}

\begin{theorem}

\label{DedekindNP}

Let $(G,+)$ be an Abelian group and $D$ a Dedekind domain. Then the following statements hold:

\begin{enumerate}

\item The semiring $(\Id(D) , + , \cdot)$ is a simple semidomain.

\item If $n\geq 2$ and $\mathfrak{p}_1 , \ldots, \mathfrak{p}_n$ are distinct maximal ideals of $D$, then $\mathfrak{p}_1 + \mathfrak{p}_2 \cdots \mathfrak{p}_n = D$.

\item A modular function $f: \Id(D) \longrightarrow G$ can be characterized by the values of $f$ at $(0)$, $D$, and all maximal ideals of $D$.

\end{enumerate}

\begin{proof}
Let $D$ be a Dedekind domain and $\mathfrak{a}$ a nonzero and proper ideal of $D$. Then $\mathfrak{a} = \mathfrak{p}_1^{\alpha_1} \mathfrak{p}_2^{\alpha_2} \cdots \mathfrak{p}_k^{\alpha_k}$, where $\mathfrak{p}_i$ are maximal ideals of $D$, and $\alpha_i >0, k \geq 2$ are integers. Then by Lemma \ref{simplesemiring}, we have \[f(\mathfrak{a})=f(\mathfrak{p}_1^{\alpha_1} \mathfrak{p}_2^{\alpha_2} \cdots \mathfrak{p}_k^{\alpha_k})=f(\mathfrak{p}_1 \mathfrak{p}_2 \cdots \mathfrak{p}_k).\]
	Now since $\mathfrak{p}_i$s are maximal, modularity of $f$ gives the following equality:
		\begin{align*}
			f(D)+f(\mathfrak{p}_1 \mathfrak{p}_2 \cdots \mathfrak{p}_k) &= f(\mathfrak{p}_1)+f(\mathfrak{p}_2 \mathfrak{p}_3 \cdots \mathfrak{p}_k).
		\end{align*}
	Applying modular property for the ideals $\mathfrak{p}_2$ and $\mathfrak{p}_3 \mathfrak{p}_4 \cdots \mathfrak{p}_k$, we get
		\begin{align*}
			2f(D)+f(\mathfrak{p}_1 \mathfrak{p}_2 \cdots \mathfrak{p}_k) &= f(\mathfrak{p}_1) + \left(f(D)+f(\mathfrak{p}_2 \mathfrak{p}_3 \cdots \mathfrak{p}_k)\right)\\
			& = f(\mathfrak{p}_1)+f(\mathfrak{p}_2)+f(\mathfrak{p}_3 \cdots \mathfrak{p}_k).
		\end{align*}
	Continuing this process, we can write
		\begin{align}\label{eq3}
			f(\mathfrak{a})=f(\mathfrak{p}_1 \mathfrak{p}_2 \cdots \mathfrak{p}_k) = f(\mathfrak{p}_1)+f(\mathfrak{p}_2) + \cdots + f(\mathfrak{p}_k)-(k-1)f(D).
		\end{align}

Our claim is that if the values of $f$ at $(0)$, $D$, and all maximal ideals of $D$ are given, we can define $f$ for an arbitrary nonzero proper ideal $\mathfrak{a} = \mathfrak{p}_1^{\alpha_1} \mathfrak{p}_2^{\alpha_2} \cdots \mathfrak{p}_k^{\alpha_k}$ ($\alpha_i >0$ and $k \geq 1$) of $D$ as follows:

\begin{align}
			f(\mathfrak{a}) = f(\mathfrak{p}_1)+f(\mathfrak{p}_2) + \cdots + f(\mathfrak{p}_k)-(k-1)f(D).
		\end{align}

and then $f$ becomes modular.

In order to show the trueness of our claim, we distinguish four cases for the following formula:

\begin{align}\label{eq4}
f(\mathfrak{a} + \mathfrak{b}) + f(\mathfrak{a} \cdot \mathfrak{b}) = f(\mathfrak{a}) + f(\mathfrak{b}).
\end{align}

\begin{enumerate}

\item If one of the ideals $\mathfrak{a}$ and $\mathfrak{b}$ is the zero ideal, then it is straightforward that equality (\ref{eq4}) holds.

\item Also, if one of the ideals $\mathfrak{a}$ and $\mathfrak{b}$ is the whole domain $D$, again it is straightforward that equality (\ref{eq4}) holds.

\item By Proposition 6.8 in \cite{Karpilovsky1989}, if the ideals $\mathfrak{a} = \mathfrak{p}_1^{\alpha_1} \mathfrak{p}_2^{\alpha_2} \cdots \mathfrak{p}_k^{\alpha_k}$ and $\mathfrak{b} = \mathfrak{q}_1^{\beta_1} \mathfrak{q}_2^{\beta_2} \cdots \mathfrak{q}_l^{\beta_l}$ have no common maximal factor in their decompositions, then the ideals $\mathfrak{a}$ and $\mathfrak{b}$ are comaximal, i.e., $$\mathfrak{a} + \mathfrak{b} = D,$$ and we have the following equalities: $$f(\mathfrak{a} + \mathfrak{b}) = f(D),$$  $$f(\mathfrak{a}) = f(\mathfrak{p}_1)+f(\mathfrak{p}_2) + \cdots + f(\mathfrak{p}_k)-(k-1)f(D),$$ $$f(\mathfrak{b}) = f(\mathfrak{q}_1)+f(\mathfrak{q}_2) + \cdots + f(\mathfrak{q}_l)-(l-1)f(D),$$ $$f(\mathfrak{ab}) = f(\mathfrak{p}_1)+f(\mathfrak{p}_2) + \cdots + f(\mathfrak{p}_k) + f(\mathfrak{q}_1)+f(\mathfrak{q}_2) + \cdots + f(\mathfrak{q}_l) - (k+l-1)f(D), $$ which show that equality (\ref{eq4}) holds in this case as well.

\item By Proposition 6.8 in \cite{Karpilovsky1989}, if the ideals $\mathfrak{a} = \mathfrak{p}_1^{\alpha_1} \mathfrak{p}_2^{\alpha_2} \cdots \mathfrak{p}_k^{\alpha_k}$ and $\mathfrak{b} = \mathfrak{q}_1^{\beta_1} \mathfrak{q}_2^{\beta_2} \cdots \mathfrak{q}_l^{\beta_l}$ have common maximal factors in their decompositions and those common maximal factors are the ideals $\mathfrak{r}_1, \mathfrak{r}_2, \ldots, \mathfrak{r}_t$, then $\mathfrak{a} + \mathfrak{b} = \mathfrak{r}_1^{\gamma_1} \mathfrak{r}_2^{\gamma_2} \cdots \mathfrak{r}_t^{\gamma_t}$ and the following equalities satisfy: $$f(\mathfrak{a} + \mathfrak{b}) = f(\mathfrak{r}_1)+f(\mathfrak{r}_2) + \cdots + f(\mathfrak{r}_t)-(t-1)f(D),$$  $$f(\mathfrak{a}) = f(\mathfrak{p}_1)+f(\mathfrak{p}_2) + \cdots + f(\mathfrak{p}_k)-(k-1)f(D),$$ $$f(\mathfrak{b}) = f(\mathfrak{q}_1)+f(\mathfrak{q}_2) + \cdots + f(\mathfrak{q}_l)-(l-1)f(D),$$ $$f(\mathfrak{ab}) = [\sum_{i=1}^k f(\mathfrak{p}_i)+ \sum_{j=1}^l f(\mathfrak{q}_j) - \sum_{s=1}^t f(\mathfrak{r}_t)] - (k+l-t-1)f(D), $$
\end{enumerate}
and the proof is complete.
\end{proof}

\end{theorem}

\begin{corollary}

\label{parvardinaturalnumbers}

Let $(G,+)$ be an Abelian group and $f$ be a function from $\mathbb N_0$ into $G$ with the following property: $$f(\lcm\{a,b\}) + f(ab) = f(a) + f(b).$$ Then $f$ can be characterized by the values of $f$ at numbers $0$, $1$, and all prime numbers of $\mathbb N$.

\begin{proof}
The semiring $(\mathbb N_0, \lcm, \cdot)$ is isomorphic to the semiring $(\Id(\mathbb Z), +, \cdot)$ and $\mathbb Z$ is a Dedekind domain.
\end{proof}
\end{corollary}

\begin{example}

\label{DedekindNPExample}

Let $D$ be a Dedekind domain. Define a function $f: D \longrightarrow \mathbb Z$ as follows:

\begin{enumerate}

\item $f((0)) = f(D) = 0$.

\item If $\mathfrak{a} = \mathfrak{m}_1^{\alpha_1} \cdots \mathfrak{m}_k^{\alpha_k}$, where $\mathfrak{m}_1, \ldots, \mathfrak{m}_k$ are distinct maximal ideals of $D$, then $f(\mathfrak{a}) = k$.

\end{enumerate}

By Theorem \ref{DedekindNP}, $f$ is modular and one can interpret $f$ as a function that counts the number of maximal factors in the decomposition of an arbitrary ideal of $D$.

\end{example}

Now we go further to give other important examples for nontrivial modular functions in different branches of mathematics:

\begin{example}

\label{modularexamples}

Some examples for modular functions in different branches of mathematics:

\begin{enumerate}

\item Let $\Omega$ be a finite set and $2^{\Omega}$ the set of all subsets of $\Omega$. The size function $|\cdot| : 2^{\Omega} \longrightarrow \mathbb N_0$ is obviously a modular function. On the other hand, according to classical probability, the probability of an event $A \subseteq \Omega$ is defined by $P(A) = |A| / |\Omega|$. So, $P(A \cup B) + P(A \cap B) = P(A) + P(B)$ for any $A,B \subseteq \Omega$. This means that $P: 2^{\Omega} \longrightarrow \mathbb [0, +\infty)$ is also a modular function \cite{Stirzaker2003}. It is clear that $(2^{\Omega}, \cup, \cap, \emptyset, \Omega, ^{\prime})$ is a Boolean algebra and therefore a semiring.

\item Let $V$ be a finite dimensional inner-product vector space. If we denote the set of all subspaces of $V$ by $\Sub(V)$, then the function $\dim : \Sub(V) \longrightarrow \mathbb N_0$ is a modular function, since by Proposition 5.16 in \cite{Golan2007}, we have $$\dim(W_1 + W_2) + \dim(W_1 \cap W_2) = \dim(W_1) + \dim(W_2).$$ It is clear that $(\Sub(V), +, \cap, (0), V, ^{\bot})$ is a Boolean algebra and therefore a semiring.

\item Let $M$ be an $R$-module of finite length and let $K$ and $N$ be $R$-submodules of $M$. Then $c(K+N) + c(K \cap N) = c(K) + c(N)$, where $c(L)$ is the composition length of a module \cite[Corollary 11.5]{AndersonFuller1992}. Now, if we consider $M$ to be a distributive module \cite{Camillo1975} of finite length, then $\Sub(M)$ is a bounded distributive lattice and therefore, a semiring. For an example of a distributive module of finite length, refer to \cite{DlabRingel1972}.

\item Let $D$ be an integral domain. Then $D$ is a Pr\"{u}fer domain if and only if $(I + J)(I \cap J) = IJ$ for all ideals $I,J$ of $D$. On the other hand, if $D$ is a Pr\"{u}fer domain, then $(\Id(D), + , \cap)$ is a semiring \cite[Theorem 6.6]{LarsenMcCarthy1971}, where by $\Id(D)$, we mean the set of all ideals of $D$. Now it is clear that the function $f$ from the semiring $(\Id(D), + , \cap)$ to the commutative semigroup $(\Id(D), \cdot)$ defined by $f(I) = I$ is a modular function.

\item Let $(G,+)$ be a partially-ordered Abelian group. It is easy to see that $a \vee b$ exists if and only if $a \wedge b$ exists in $G$ and in this case, by Theorem 15.2 in \cite{Gilmer1972}, $$a \vee b + a \wedge b = a + b.$$ Therefore, in such a case, the function $f: (G,\vee,\wedge) \longrightarrow (G,+)$ defined by $f(x) = x$ is a modular function. Now, if $(G,\vee,\wedge)$ is also a distributive lattice, $f$ is a modular function over the pre-semiring $G$.

\item Let $L$ be a bounded distributive lattice of finite length. Then the function $\height: L \longrightarrow \mathbb N_0$ is a modular function, since $$\height(a\vee b) + \height(a\wedge b) = \height(a) + \height(b)$$ for all $a,b \in L$ \cite[Corollary 376]{Graetzer2011}.

\item Consider the Peano-Jordan content $c$ (volume, area or length) over an Euclidian space $\mathbb E^n$. Define $FC(\mathbb E^n)$ to be the set of all subsets $A \subseteq \mathbb E^n$ such that $c(A) < \infty$. Then $(FC(\mathbb E^n), \cup, \cap, \emptyset)$ is a hemiring and $c$ is a modular function over $FC(\mathbb E^n)$ \cite[Theorem 22]{Rogosinski1952}.

\end{enumerate}

\end{example}

Now, we give a pre-semiring version of the so-called inclusion-exclusion theorem in probability theory. Note that the set-theoretic version of the following theorem is due to Poincar\'{e} \cite[Exercise 1.12.52.]{Bogachev2007}. It is also a good point to mention that one can see the number theory version of the inclusion-exclusion theorem in Theorems 7.3 and 7.4 in \cite{Loo-Keng1982}. Also, for a multiplicative ideal theory version of this theorem refer to the recent paper \cite{Anderson2016} by D. D. Anderson et al. on GCD and LCM-like identities for ideals in commutative rings.

\begin{theorem}[Inclusion-Exclusion Theorem]

\label{Poincare}

Let $S$ be a multiplicatively idempotent pre-semiring and $T$ a commutative semigroup. If $m: S \longrightarrow T$ is a modular function, then for any elements $s_1, s_2,\ldots, s_n$ in S, the following equality holds: $$m\left(\sum^n_{i=1} s_i\right) + \sum_{i<j} m(s_i s_j) + \sum_{i<j<k<l} m(s_i s_j s_k s_l) + \cdots =$$ $$\sum^n_{i=1} m(s_i) + \sum_{i<j<k} m(s_i s_j s_k)+ \cdots \qquad \text{(Poincar\'{e}'s Formula)},$$ where the number of multiplicative factors in the sums of both sides of the equality is at most $n$.

\begin{proof} The proof is by mathematical induction. Since $m$ is modular, the Poincar\'{e}'s Formula holds for $k=2$. Now let Poincar\'{e}'s Formula hold for $k=n$ and take $s_1, \ldots, s_n, s_{n+1} \in S$. Since $m$ is modular, we have $$m\big(\sum^n_{i=1} s_i + s_{n+1}\big) + m\big((\sum^n_{i=1} s_i) \cdot s_{n+1}\big) = \sum^n_{i=1} m(s_i) + m(s_{n+1}).$$

Now by applying Poincar\'{e}'s Formula for the $n$ elements $s_1 s_{n+1}, s_2 s_{n+1},\ldots, s_n s_{n+1}$, we have $$m\big(\sum^n_{i=1} s_i s_{n+1}\big) + \sum_{i<j} m(s_i s_{n+1} s_j s_{n+1}) + \sum_{i<j<k<l} m(s_i s_{n+1} s_j s_{n+1} s_k s_{n+1} s_l s_{n+1}) + \cdots =$$ $$\sum^n_{i=1} m(s_i s_{n+1}) + \sum_{i<j<k} m(s_i s_{n+1} s_j s_{n+1} s_k s_{n+1})+ \cdots .$$ By using that $S$ is multiplicatively idempotent, the above equality can be written as follows:

$$m\big(\sum^n_{i=1} s_i s_{n+1}\big) + \sum_{i<j} m(s_i s_j s_{n+1}) + \sum_{i<j<k<l} m(s_i s_j s_k s_l s_{n+1}) + \cdots =$$ $$\sum^n_{i=1} m(s_i s_{n+1}) + \sum_{i<j<k} m(s_i s_j s_k s_{n+1})+ \cdots .$$ Finally, by adding $m(\sum^n_{i=1} s_i + s_{n+1})$ to both sides of the recent equality, we get the Poincar\'{e}'s Formula for the $n+1$ elements $s_1, s_2, \ldots, s_{n+1}$ and the proof is complete.
\end{proof}

\end{theorem}

\begin{proposition}

\label{Poincare2}

Let $S$ be a multiplicatively idempotent pre-semiring, $R$ a ring, and $f: S \longrightarrow R$ a modular function. For any $n>1$, if $s_1,s_2,\ldots,s_n$ in $S$ are independent, then $s_1 + s_2 + \cdots + s_{n-1}$ and $s_n$ are also independent.

\begin{proof}

 Since $S$ is multiplicatively idempotent, by Theorem \ref{Poincare}, we have: $$f\big((s_1 + \cdots + s_{n-1})s_n\big) = f(s_1 s_n + \cdots + s_{n-1} s_n) = $$ $$\sum^{n-1}_{i=1} f(s_i s_n) - \sum_{1\leq i<j <n} f(s_i s_j s_n) + \sum_{1\leq i<j<k <n} f(s_i s_j s_k s_n) - \cdots + (-1)^n f(s_1 s_2 \cdots s_n) =$$ $$\big[\sum^{n-1}_{i=1} f(s_i) - \sum_{1\leq i<j <n} f(s_i s_j) + \sum_{1\leq i<j<k <n} f(s_i s_j s_k) - \cdots + (-1)^n f(s_1 s_2 \cdots s_{n-1})\big] f(s_n) =$$
$$f(s_1 + \cdots + s_{n-1}) f(s_n),$$ and the proof is complete.
\end{proof}

\end{proposition}

In the following theorem and its corollaries, we show when a finitely additive function becomes modular:

\begin{theorem}

\label{modular}

Let $S$ be a zerosumfree semiring such that $1+1 \in \comp(S)$, $T$ a commutative semigroup, and $f: S \longrightarrow T$ a finitely additive function. Then $f\mid_{\comp(S)} : \comp(S) \longrightarrow T$ is modular.

\begin{proof}
Let $s,t \in \comp(S)$. Therefore by Proposition \ref{complemented1} and Theorem \ref{zerosumfree}, we have $f(s+t) = f(s \sqcup t) = f(s + s^{\bot} t) = f(s) + f(s^{\bot} t)$. Finally, by Proposition \ref{complemented1}, $f(s+t) + f(st) = f(s) + f(s^{\bot} t) + f(st) = f(s) + f(t)$ and the proof is complete.
\end{proof}

\end{theorem}

By Theorem \ref{zerosumfree}, it is clear that if each element of a semiring $S$ is complemented, then $S$ is a Boolean algebra. Therefore, it is obvious that we have the following result:

\begin{corollary}

Let $S$ be a Boolean algebra and $T$ a commutative and cancellative monoid. Then the following statements for a function $f: S \longrightarrow T$ are equivalent:

\begin{enumerate}

\item $f$ is finitely additive,

\item $f$ is modular and $f(0) = 0$.

\end{enumerate}

\end{corollary}

Note that if $L$ is a Stone algebra, then $\Skel(L)$ is a Boolean algebra \cite[Theorem 214]{Graetzer2011}. Therefore, we have the following:

\begin{corollary}
Let $L$ be a Stone algebra and $T$ a commutative and cancelation monoid. If $f: L \longrightarrow T$ is finitely additive, then $f\mid_{\Skel(S)} : \Skel(S) \longrightarrow T$ is modular.
\end{corollary}

\begin{proposition}

\label{Poincare3}

Let $S$ be a zerosumfree semiring such that $1+1 \in \comp(S)$, $T$ a commutative semigroup, and $f: S \longrightarrow T$ a finitely additive function. For any $s_1, s_2,\ldots, s_n \in \comp(S)$, we have the following: $$m\big(\sum^n_{i=1} s_i\big) + \sum_{i<j} m(s_i s_j) + \sum_{i<j<k<l} m(s_i s_j s_k s_l) + \cdots =$$ $$\sum^n_{i=1} m(s_i) + \sum_{i<j<k} m(s_i s_j s_k)+ \cdots \qquad \text{(Poincar\'{e}'s Formula)},$$ where the number of multiplicative factors in the sums of both sides of the equality is at most $n$.

\begin{proof}

Let $S$ be a zerosumfree semiring such that $1+1 \in \comp(S)$. So by Theorem \ref{zerosumfree}, $\comp(S)$ is a Boolean algebra. This means that $\comp(S)$ is a multiplicatively idempotent semiring. On the other hand, by Theorem \ref{modular}, $m\mid_{\comp(S)} : \comp(S) \longrightarrow T$ is a modular function. Hence, by Theorem \ref{Poincare}, Poincar\'{e}'s Formula holds for complemented elements of $S$ and this completes the proof.
\end{proof}

\end{proposition}

\section*{Acknowledgment} The first named author is supported by Department of Engineering Science at Golpayegan University of Technology. The authors wish to thank Sepanta Asadi for informing them about a point related to Proposition \ref{constantmodular2}. They are also grateful for the useful comments of the anonymous referee.

\end{document}